\documentclass[12pt]{article}
\usepackage{pstricks}
\usepackage{graphicx}%,psfrag}
\usepackage{amsmath}
\usepackage{amssymb}
\usepackage{amsthm}
\makeatletter

\renewcommand{\@makecaption}[2]{%
\vspace{\abovecaptionskip}%
\sbox{\@tempboxa}{#1. #2}
#1. #2\par
\vspace{\belowcaptionskip}} %Left flushing for figures' captions.
\makeatother
\begin{document}

\begin{center}
\textbf{\Large Motion of a rough disc in Newtonian aerodynamics}

\bigskip

{\large Sergey Kryzhevich${}^{a,b}$\footnote{Email: kryzhevicz@gmail.com, s.kryzhevich@spbu.ru}}% and Alexander Plakhov${}^{a,d}$}

\textit{\noindent\small ${}^a$ Center for Research and Development in Mathematics and Applications, Department of Mathematics, University of Aveiro, Aveiro 3810-193, Portugal;\\
${}^b$ Faculty of Mathematics and Mechanics, Saint-Petersburg State University, 28, Universitetskiy pr., Peterhof, Saint-Petersburg, 198503, Russia.}

\end{center}

\bigskip

 %\begin{center}

{\small \textbf{Abstract.} Dynamics of a rough disc in a rarefied medium is considered.  The main result of the paper is the following: any finite rectifiable curve can be approximated in the Hausdorff metric by trajectories of centers of rough discs (that is, $C^0$-small perturbations of regular discs), provided that the parameters of the system are carefully chosen. To control the dynamics of the disc, we use the so-called inverse Magnus effect which causes deviation of the trajectory of a spinning body. We use the shape of the perturbed disc as a parameter to control the Magnus effect. We study the so-called response laws for scattering billiards e.g. relationship between the velocity of incidence and that of reflection. We construct a special family of such laws that is weakly dense in the set of symmetric Borel measures. Then we find a shape of cavities that provides selected law of reflections.  Mathematical models for dynamics of rough bodies, corresponding to cavities of the constructed shape, are studied. We write down differential equations that describe motions of such discs. We offer a method on how a given curve can be approximated by considered trajectories.
§§
\bigskip

\noindent\textbf{Keywords:} Billiards; shape optimization; Magnus effect; rarefied medium; retroreflectors.}

\bigskip

{\small \noindent{\textbf{Mathematics subject classifications:} 37D50, 49Q10, 70Q05, 93B05.}}

%\end{center}

\section*{Acknowledgements}

\hspace{\parindent} This work was supported by Russian Foundation for Basic Researches under Grants 12-01-00275-a and 14-01-00202-a, by Saint-Petersburg State University under Thematic Plans 6.0.112.2010 and 6.38.223.2014, by FEDER funds through COMPETE~-- Operational Programme Factors of Competitiveness (Programa Operacional Factores de Competitividade) and by Portuguese funds through the Center for Research and Development in Mathematics and Applications (CIDMA) from the "{\it Funda\c{c}\~{a}o para a Ci\^{e}ncia e a Tecnologia}" (FCT), cofinanced by the European Community Fund FEDER/POCTI under FCT research projects (PTDC/MAT/113470/2009 and PEst-C/MAT/UI4106/2011 with COMPETE number FCOMP-01-0124-FEDER-022690). The author is grateful to Prof. Alexandre Plakhov from University of Aveiro for his ideas, remarks and corrections.

\section{Introduction}

\hspace{\parindent} Consider a body with a piecewise smooth boundary moving in a two-dimensional rarefied homogeneous medium. The particles composing this medium are initially at rest. They never interact, they collide with the body in the perfectly elastic fashion and move freely between consecutive reflections from the boundary of the body. 

This simple aerodynamic model was first introduced by Newton in his {\it Principia} (1687). He studied a particular case of this model where a convex axially symmetric body translates along its axis of symmetry. Due to collisions with particles of the medium, the force of resistance slows down the motion of the body. Newton studied the problem of finding the shape of the body that minimizes the force of resistance. The solution looks like a truncated cone with a slightly inflated lateral surface. Several generalizations of Newton's problem related to (generally) nonconvex and/or non-symmetric bodies have been studied in 1990s and 2000s by various authors [1--10]. There are open problems in this area; for instance, the shape of the convex and non-symmetric body of least resistance is not understood.

These investigations are closely related to the so-called problem of invisibility. One constructs a system of mirrors, invisible for an observer (observers), placed in a fixed point (points) or looking from a fixed direction (directions). Though the complete invisibility is impossible,  some of related problems, for example, invisibility from one point or invisibility from one direction have been already solved [11--13]. 

Dynamics of a rarefied gas is a well-studied problem, see [14--17] and references therein. 

Even more difficult and diverse are problems related to combined translational and rotational motion of bodies in a rarefied medium. Some of these problems are addressed in [20--27] under the assumption that the rotational motion is much slower than the translational one. In this case interaction of each individual particle with the body occurs as if there were no rotation at all: the turn of the body during the time of interaction can be neglected. It is shown, in particular, that the resistance of a convex body, in Euclidean space of arbitrary dimension, can be both increased and decreased by roughening its surface. The rates of maximum increase and decrease are found to depend only on the dimension and not on the original convex body; in the 3D case they are equal, respectively, to 2 and (approx.) $0.969445$.

The Newtonian dynamics of a body that performs both translational and rotational motion is a very intriguing and completely unexplored subject even in the 2D case.  Even attempts to study dynamics of very simple bodies, like a rod, not to say about an ellipse or a triangle, meet serious difficulties. The only exception is a circle, whose dynamics is trivial: the path of its center is a straight line. In this paper, we condider the so-called rough discs that represent the idea of a set, close to a ball in the Hausdorff metrics (see the beginning of Section 2). The principal goal of the article is the following.

We show that trajectories of centers of rough discs are dense in the set of finite rectifiable plane curves endowed with the Hausdorff metrics.

The proof of this statement is based on the following idea.  In the typical case, if a disc rapidly rotates, say, counterclockwise, then the velocity vector of its center of mass changes in the clockwise direction. This phenomenon is called the {\it inverse Magnus effect}, see Figure\, 1. The word "inverse"\ means this effect is inverse to the Magnus effect proper for classical gas dynamics and well-known for soccer or ping-pong players where a ball deviates at the direction of rotation. There is no contradiction: influence of a classic gas is very different from one of rarefied media.

\begin{figure}[ht]
\begin{center}
\includegraphics*[width=2.5in]{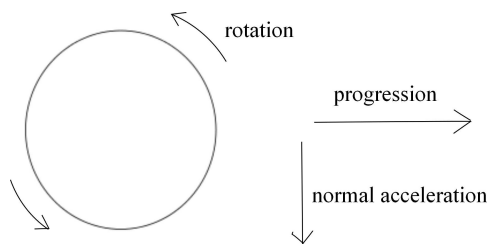}
\end{center}
\caption{Inverse Magnus effect.}
\end{figure}

The magnitude of the effect depends on the shape of cavities on the boundary of the disc and on the relative angular velocity $\lambda$ of the disc. In this paper we construct a very special cavity in such a way that (i) the relative angular velocity monotonously increases and (ii) the magnitude of the effect is nearly zero for all values of $\lambda$ except for several (relatively small) intervals of values. On these intervals the effect is adjusted so as to ensure right turn of the velocity vector to a certain angle.

Therefore, the basic idea of the proof of our main results is to use shapes of cavities to control inverse Magnus effect. We believe that our construction can be generalized to three dimensions, but postpone the 3D study to the future.

\noindent\textbf{Structure of the paper}. In Section 2, we consider an immobile scattering billiard which gives a simplified model for dynamics of a particle inside a hollow. In the next section we formulate physical assumptions on a moving body and medium and introduce some notions. We introduce the concept of $\delta$~-- pseudotrajectory, corresponding to immobile billiard system. We show that, for sufficiently precise approximations to so-called perfect rough discs, scattering billiard model gives a good approximation for relative motions of particles inside cavities. Then, in Section 4, we can apply the model for motions of rough bodies [12]. We study some special types of cavities and related reflection laws (Sections 5--7). We formulate the main result of the paper (Section 8) and prove it. We main idea of the proof is approximation of a curve by broken lines, for which we can write down equations of motions and shapes of cavities explicitly. 

\section{Laws of scattering for immobile billiards}

\hspace{\parindent} We start with the definition of a rough disc. Fix $r > 0$, take a regular $n_0$-gon $(n_0 \ge 3)$ inscribed in a circle with radius $r$ (let its center be $O$), and replace each side of the $n_0$-gon with a curve joining its endpoints. Each curve is piecewise smooth, does not have any self-intersections, and is contained in the circular sector with vertices at $O$ and at the endpoints of the corresponding side. In addition, all the curves are congruent: each curve can be obtained from another one by rotation around $O$ by $2\pi k/n_0$. For each integer $n > n_0$ make a similar procedure: take a regular $n$-gon inscribed in the same circle and replace its sides with curves, so that the obtained sequence of sets tends to the circle in Hausdorff metrics. The union of the curves in each $n$-gon bounds a domain $B_n$. 

\noindent\textbf{Definition 1.} The sequence of domains $B_n,\, n \ge n_0$ is called a {\it rough disc}.

Thus, a rough disc is an idealized object. It is not a domain, but rather it can be informally viewed as the "limit" of a sequence of domains $B_n$. Its "boundary" is obtained by repetition of identical infinitesimal curves similar to the original one. They are interpreted as infinitesimal hollows on the disc boundary. The billiard scattering by the rough disc is uniquely defined by the shape of the curve.

We assume that the dics moves in a medium where the mass is uniformly distributed according to the measure $m$:  $dm= \rho \, d{\mathrm S}$, where ${\mathrm S}$ is the Lebesgue measure in ${\mathbb R}^2$. We treat particles as infinitesimal parts of the medium. We neglect Brownian motion of particles when we calculate interactions between particles and the body. However, we suppose that particles instantly fill in the space after the body was passing.  Let $\boldsymbol{X}=(X,Y)$ be the current position of the center of the body, $\phi$ be the current angle of rotation of the body with respect to its initial position. Let $\boldsymbol V$ be the velocity of the center, $|{\boldsymbol V}|=V$. Denote by $\omega$ the angular velocity, by $r$ be the radius of the disc and by $I=\kappa M r^2$ the moment of inertia. For a regular disc with uniformly distributed mass we have $I=Mr^2/2$, and in any case $\kappa\in [0,1]$. We introduce the angular coordinate $\xi$ on the boundary of the ordinary disc representing the smooth approximation of the moving body, identifying this boundary with the unit circle $S^1=[-\pi,\pi]/\{\pi=-\pi\}$. Recall the notion for the dimensionless relative angular velocity $\lambda=\omega r/V$. The force of resistance of the medium acting on the disc and the moment of this force are defined as limits, when $n \to \infty$, of the force and the moment of force acting on $B_n$. Using these values, we derive the equations of motion of a rough disc on the plane. These equations, and therefore the trajectory of the disc, depend on the shape of the infinitesimal curve forming its boundary. The natural question arises: which curves can be traversed by the disc center?

Description of scattering by a rough disc and equations of motions for such discs can be found in [27], and in chapters 4 and 7 of the book [16]. We partly reproduce them here.

\noindent\textbf{Definition 2.} A {\it hollow} is a piecewise smooth non self-intersecting curve contained in a closed isosceles triangle whose base is the segment joining the endpoints of the curve. The segment is called the {\it opening} of the hollow.

We use the notion $\Omega$ for a hollow and $I$ for its opening. Introduce the uniform coordinate $\xi \in [0,\, 1]$ on the opening $I$; the values $\xi = 0$ and $\xi = 1$ correspond to its endpoints. Let $\boldsymbol n$ be the unit outer normal to $I$. Consider a particle that enters a hollow $\Omega$ through its opening $I$. Fix the point $\xi$ where it intersects the opening and the angle $\varphi\in (-\pi/2,\pi/2)$ formed by $-\boldsymbol{n}$ and the incidence velocity $\boldsymbol{v}$. If the particle makes a finite number of reflections from regular points of $\Omega$, intersects $I$ again and leaves, we denote by $\xi^+ = \xi^+_\Omega(\varphi,\xi)$ the point of the second intersection and by $\varphi^+ = \varphi^+_\Omega(\varphi,\xi)$ the angle formed by $\boldsymbol{n}$ and the velocity ${\boldsymbol v}^+$. 

Almost all particles leave the hollow $\Omega$ after a finite number of reflections. This follows from the measure-preserving property of billiard and from  Poincar\'e's recurrence theorem. Thus for almost all initial conditions $(\varphi,\xi) \in [-\pi/2,\, \pi/2] \times [0,\, 1]$ the values $\varphi^+_\Omega(\varphi,\xi)$ and $\xi^+_\Omega(\varphi,\xi)$ are well-defined. Introduce the probability measure $\mu$ on $[-\pi/2,\, \pi/2] \times [0,\, 1]$ according to $d\mu(\varphi,\xi) = \frac 12\, \cos\varphi\, d\varphi\, d\xi$. The map $T_\Omega : (\varphi,\xi) \mapsto (\varphi^+_\Omega(\varphi,\xi),\, \xi^+_\Omega(\varphi,\xi))$ is defined on a full-measure subset of $[-\pi/2,\, \pi/2] \times [0,\, 1]$ and maps it bijectively onto itself. Moreover, it preserves the measure $\mu$ and is involutive, $T_\Omega = T_\Omega^{-1}$.

Next introduce the Borel measure $\eta_\Omega$ on the square $\Box := [-\pi/2,\, \pi/2] \times [-\pi/2,\, \pi/2]$ as follows: $\eta_\Omega(A) = \mu(\{ (\varphi,\xi) : (\varphi,\, \varphi^+_\Omega(\varphi,\xi)) \in A \})$ for any Borel set
$A \subset \Box$. 

This measure can be defined in a different way: let $\sigma_\Omega$ be the mapping $(\varphi,\xi) \mapsto (\varphi,\, \varphi^+_\Omega(\varphi,\xi))$
from $[-\pi/2,\, \pi/2] \times [0,\, 1]$ to $\Box$; then $\eta_\Omega$ is the measure $\eta_\Omega = \sigma_\Omega^\# \mu$. Here $\sigma_\Omega^\#$ is the mapping of measures induced by $\sigma_\Omega$.

\noindent\textbf{Definition 3.} $\eta_\Omega$ is called the {\it measure induced by the hollow} $\Omega$.

Define the probability measure $\gamma$ on $[-\pi/2,\, \pi/2]$ by $d\gamma(\varphi) = \frac 12\, \cos\varphi\, d\varphi$. For a set $A \subset \Box$, denote $A^* = \{ (\varphi,\varphi^+) : (\varphi^+,\varphi) \in A \}$. Denote by $\Upsilon$ the set of Borel measures $\eta$ on $\Box$ such that for  all Borel sets $A \subset \Box$ and $I \subset [-\pi/2,\, \pi/2]$ one has $\eta(A) = \eta(A^*)$  and $\eta(I \times [-\pi/2,\, \pi/2]) = \gamma(I)$. The fact that $\eta_\Omega \in \Upsilon$ can be easily deduced from the measure preserving and involutive properties of the map $T_\Omega$; see [12] for details. The following important theorem states that, inversely, the set of measures induced by hollows is weakly dense in $\Upsilon$.

\noindent\textbf{Density Theorem [12].} \emph{The set $\{ \eta_\Omega : \Omega \text{ is a hollow\,} \}$ is weakly dense in $\Upsilon$. In other words, for any $\eta \in \Upsilon$ there exists a sequence of hollows $\Omega_k$ such that
$$\lim_{k\to\infty} \int\!\!\!\int_\Box f(\varphi,\varphi^+)\, d\eta_{\Omega_k}(\varphi,\varphi^+) = \int\!\!\!\int_\Box f(\varphi,\varphi^+)\, d\eta(\varphi,\varphi^+)$$
for any continuous function $f : \Box \to {\mathbb R}$.}

\section{Pseudotrajectories}

\hspace{\parindent} Given a $\delta>0$, we introduce the concept of a $\delta$-\emph{pseudotrajectory} for a billiard. 

\noindent\textbf{Definition 4.} We say that a piecewise $C^1$ smooth curve ${\boldsymbol x}(t):$ $t\in [t_0, {\hat t}_0]$ is a $\delta$-pseudotrajectory for the exterior billiard corresponding to an immobile body $A$ if the following statements are true.
\begin{enumerate}
\item ${\boldsymbol x}(t)\notin \mbox{\rm int}\, A$ for all $t\in [t_0, {\hat t}_0]$.
\item The set of $t$ such that ${\boldsymbol x}(t)\in \partial A$ is finite. Let it be $\{t_1,\ldots,t_N\}: N\in {\mathbb N}\bigcup \{0\}$. We also use the notation $t_{N+1}={\hat t}_0$.
\item For all $k\in\{1,\ldots,N\}$ the velocities ${\boldsymbol v}_{r+} = {\boldsymbol v}(t_k+0) = \dot{{\boldsymbol x}}(t_k+0)$ and ${\boldsymbol v}_{r-} = {\boldsymbol v}(t_k-0) = \dot{{\boldsymbol x}}(t_k-0)$ of the corresponding impacts satisfy inequalities
$$|{\boldsymbol v}_{e+} - {\boldsymbol v}_{e-} - 2\langle {\boldsymbol v}_{e-},\, {\boldsymbol n} \rangle {\boldsymbol n}|\le \delta. \eqno (1)$$
If $x(t)$ is a singularity point, we select one of two possible values for normal vectors.
\item The function ${\boldsymbol v}(t)=\Dot{\boldsymbol x}(t)$ is piecewise smooth and $|{\boldsymbol v}(t)-{\boldsymbol v}(t_k+0)|\le \delta$ for any $k\in \{0,\ldots,N\}$ and any $t\in (t_k,t_{k+1})$.
\end{enumerate}

We use this notion to describe trajectories of particles of non-zero mass that interact with a moving and rotating body.

\noindent\textbf{Definition 5.} A rough disc defined by a sequence $B_n$ is \emph{perfect} if there exist $m_0\in {\mathbb N}$, $\lambda_0>0$, $0<\delta<\pi/2$ and $K>0$ such that for any $m\ge m_0$ all $\delta$-pseudotrajectories, entering the corresponding hollow with the incident angle $\ge \lambda_0$, have at most $K$ impacts before they leave the hollow.

 We make the following assumptions on interactions between the body and particles.
\begin{enumerate}
\item If a particle collides once with a point of the boundary of the body out of any hollow, we neglect all later interactions between the particle and the body.
\item We assume that there is a number $K>0$ such that all particles interacting with a fixed hollow of the body, have at most $K$ impacts and leave the hollow. 
\end{enumerate}

Direct calculations lead us to the following statement.

\noindent \textbf{Lemma 1. }\emph{Let $\{B_n\}$ be a rough disc, $h_n$ be diameters of corresponding cavities. Let $a_0$, $v_0$, $v_1$ and $\omega_0$ be positive constants. Suppose that a body $B_n$ translates and rotates during a period $[0,T]$ so that  $|\Dot{\boldsymbol X}(t)|\in [v_0,v_1]$, $|\Ddot{\boldsymbol X}(t)|\le a_0$. Here ${\boldsymbol X}(t)$ is the position of the center of the body, $\omega$ is the angular velocity. Then there exist $n_0\in {\mathbb N}$ and $C_1,C_2>0$ such that any particle entering a cavity with incidence angle $\le \lambda_0$ spends at most $C_1 h_n$ units of time inside the cavity. The part of the trajectory inside a cavity forms a $C_2h_n$ pseudotrajectory.}

\section{Dynamics of perfect rough bodies}

\hspace{\parindent} The dynamics of the rough disc is described by the following system of ordinary differential equations [12, Theorem 7.1, P.\,203]:
$$
M \Dot{\boldsymbol V}={\boldsymbol R}(\eta,\omega,{\boldsymbol V})=\dfrac83 r\rho V^2 \overline{\boldsymbol R}(\eta,\lambda);\qquad
I \dot \omega=R_I(\eta,\omega,{\boldsymbol V})=\dfrac83 r\rho V^2 \overline{R_I}(\eta,\lambda). \eqno (2)$$
Here $\eta$ is the billiard law corresponding to the selected rough disc. Formulae for dimensionless resistances $\overline{\boldsymbol R}$ depend on $\lambda$. Here we assume that  $\lambda>1$. Consider the coordinate system associated with the vector ${\boldsymbol V}$ and the orthogonal vector ${\boldsymbol V}^\bot$. Functions $\overline{\boldsymbol R}$ and $R_I$ can be found from the following formulae:
$$\begin{array}{cc}
\overline{\boldsymbol R}(\eta,\lambda)=(R_T(\eta,\lambda),R_L(\eta,\lambda));&\quad
R_T(\eta,\lambda)=\int\limits_\Box c_T(x,y,\lambda)\, d\eta(x,y);\\
R_L(\eta,\lambda)=\int\limits_\Box c_L(x,y,\lambda)\, d\eta(x,y);&\quad
R_I(\eta,\lambda)=\int\limits_\Box c_I(x,y,\lambda)\, d\eta(x,y).
\end{array}
\eqno (3)$$
Here 
$$\begin{array}{c}
c_T(x,y,\lambda)=\dfrac{3\cos \frac{x-y}{2}}{\sin \zeta}((\lambda^3\sin^3 x+3\lambda\sin x\sin^2\zeta)\cos\zeta\cos\frac{x-y}{2}-\\
(3\lambda^2\sin^2 x\sin \zeta+\sin^3\zeta)\sin\zeta\sin\frac{x-y}{2})\chi_{x\ge x_0}(x,y);\\
c_L(x,y,\lambda)=-\dfrac{3\cos \frac{x-y}{2}}{\sin \zeta}((\lambda^3\sin^3 x+3\lambda\sin x\sin^2\zeta)\cos\zeta\sin\frac{x-y}{2}+\\
(3\lambda^2\sin^2 x\sin \zeta+\sin^3\zeta)\sin\zeta\cos\frac{x-y}{2})\chi_{x\ge x_0}(x,y);\\
c_I(x,y,\lambda)=-\dfrac32\dfrac{\lambda^3\sin^3 x+3\lambda\sin x\sin^2\zeta}{\sin \zeta}(\sin x+\sin y)\chi_{x\ge x_0}(x,y);
\end{array}\eqno (4)$$
$\zeta=\arcsin \sqrt{1-\lambda^2\cos^2 x}$, $x_0=\arccos (1/\lambda)$; $\chi$ stands for the characteristic function. Applying the mentionned Theorem 7.1, we use the result of Lemma 1.

Make a transformation of variables in equations $(2)$. First of all, select $\tau$ so that:
$$d\tau=\dfrac{8r\rho V}{3M}\, dt. \eqno (5)$$
Then we define $\theta$ so that ${\boldsymbol V}=V (\cos\theta,\sin\theta)$. Let $\beta=Mr^2/I=\kappa^{-1}$ be the inverse relative moment of inertia of the rough disc. It follows from equations $(2)$ that
$$\dfrac{d\lambda}{d\tau}=\beta R_I(\lambda)-\lambda R_L(\lambda),\quad
\dfrac{dV}{d\tau}=-V R_L(\lambda),\quad
\dfrac{d\theta}{d\tau}=-R_T(\lambda).\eqno (6)$$

Observe that the variable $\tau$ is a natural parametrization of the trajectory of the center of the disc. Namely, if $S(t)$ is the path passed by the center of the disc by the moment $t$ then $dS/d\tau=3M/(8\rho r)=\mbox{const}$.

In following three sections we provide a family of specially selected roughnesses and justify that the proposed model of dynamics is applicable for rough discs with such shapes of cavities. First of all we need two types of auxiliary scattering billiards.

\section{Bunimovich mushroom}

\hspace{\parindent} Let us introduce the so-called retroreflectors. Consider a family of domains
$\Theta_h\subset {\mathbb R}^2$ ($h$ is a small positive parameter) with a piecewise smooth boundary $\partial \Theta_h$ which can be represented as a disjoint union $\partial \Theta_h=\Omega_{h}\bigcup I_{h}$
where $\Omega_{h}$ and $I_{h}$ satisfy following properties.
\begin{enumerate}
\item The arc $\Omega_{h}$ is a hollow with the opening $I_h$.
\item Consider a uniform distribution of pairs $(x_-, \nu_-)\in I_h\times [-\pi/2,\pi/2]$ that is the point and the angle of incidence.  
 Let  $\nu_+$ be the angle of the last intersection between the trajectory of a particle and the segment $I_h$, then for any $\sigma>0$ the proportion of particles such that
$|\nu_++\nu_--\pi|>\sigma$ tends to zero as $h\to 0$.
\end{enumerate}

In this paper we consider so-called "Bunimovich mushroom"\ [28,29], Figure \,2. There exist other patterns of retroreflectors [12, Chapter 9].

The pattern of the mushroom, we use in this article is the following: a domain $\Theta_h$ which is a union two domains: $\Theta_{h1}$ and $\Theta_{h2}$. The first one ("pileus" of the mushroom) is a strictly convex domain which is the upper part of an ellipse, whose principal axis is horizontal. The second part of the mushroom (call it stipe) is a $b_{12}\times b_{13}$ rectangle. Let $b_{11}$ be the length of the long axis of the ellipse $\Theta_{h1}$. We assume that the tops of the stipe coincide with the foci of the ellipse. Suppose that
$${b_{12}}/{b_{11}}= 2h; \qquad {b_{13}}/{b_{12}}= h.\eqno (7)$$
We call this value $h$ imperfectness of the mushroom.

We claim that $I_{h}=(P_LP_R)$ on Figure 2 is the entrance for the considered scattering billiard. Consequently, we suppose that $\Omega_{h}$ is the rest of the boundary of the "mushroom".

\begin{figure}[ht]\begin{center}
\includegraphics*[width=2.5in]{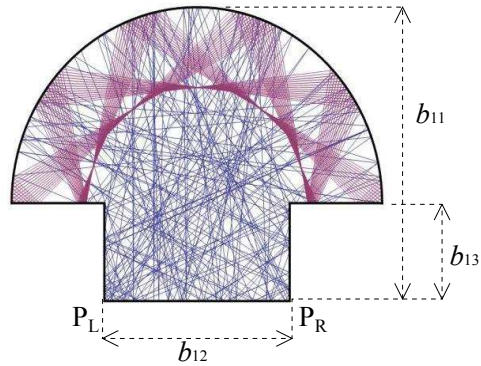}
\end{center}
\caption{Mushroom billiard.}
\end{figure}

Let us prove that this mushroom is a retroreflector.  Fixed $\sigma>0$ we define sets $\Sigma_\sigma=\{({\boldsymbol x}_-,{\boldsymbol v}_-)\in X: |{\boldsymbol x}_+-{\boldsymbol x}_-|\le \sigma b_{22},\quad |{\boldsymbol v}_++{\boldsymbol v}_-|\le \sigma\}$.

\noindent \textbf{Lemma [12, Lemma 4.1, p. 115]. }\emph{For any $\sigma>0$ there exists a $h_0>0$ such that if $h\in (0,h_0)$ and conditions $(7)$ are satisfied, the measure of the set $\Sigma_\sigma$ is greater than $1-\sigma$.}

\section{Amphora billiard: a quasi-elastic hollow}

\hspace{\parindent} Select a small positive parameter $h$ called imperfectness of the billiard. Consider two arcs of confocal parabolas  given by equations $x=\pm (1-y^2)/2$, $y\in [0,1]$. Link the lower ends of these arcs by a segment. We obtain a curvilinear triangle. Cut the middle part $F_LF_R$ of the base of this triangle corresponding to $x\in [-h,h]$. Construct two segments $A_LF_L$ and $A_RF_R$ of the length $h^2$ at ends of the obtained gap, which make angles $\pm \pi/4$ with the axis $Ox$, Figure\, 3(a)). The amphora domain is constructed.

\begin{figure}[ht]\begin{center}
\includegraphics*[width=3.6in]{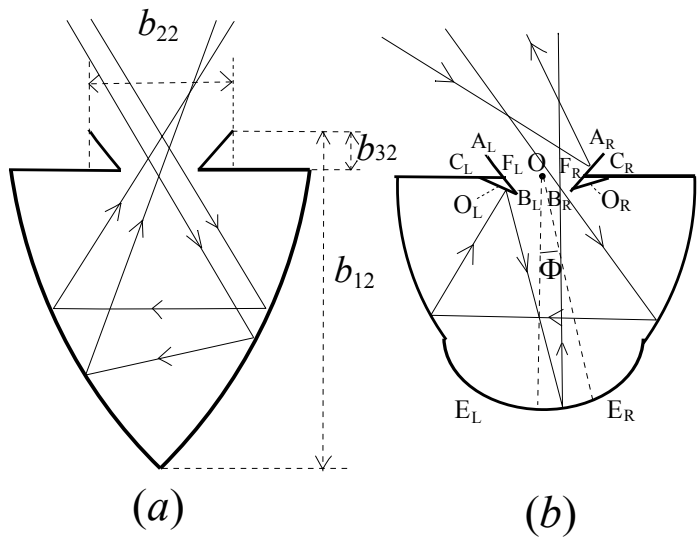}
\end{center}
\caption{Amphora billiard (a) and its modification (b).}
\end{figure}

Later on we deal with modifications of amphora billiards. We take the parameter $b_{21}$ so that
$$b_{21}= h; \qquad {b_{22}}/{b_{21}}= 2h; \qquad {b_{23}}/{b_{22}}= h \eqno (8)$$
where $h$ is the imperfectness, $b_{22}$ is the width of the entrance corridor of the billiard domain, call it "neck", $b_{23}$ is the length if this corridor.

Let $X=[-b_{22},b_{22}]\times (-\pi,0)$ be endowed with the smooth measure $\nu$ with the density $d\nu=-\sin v_-/(4b_{22})dx\,dv_-$.

Next lemma demonstrates that this amphora hollow works like a smooth mirror i.e. for "almost"\ all particles the angle of incidence "almost"\ equals to the angle of reflection. Let $N_\sigma$ be the set of initial conditions $({\boldsymbol x}_-,{\boldsymbol v}_-)\in X$ of the entrance which correspond to billiard trajectories with two impacts such that $|{\boldsymbol x}_++{\boldsymbol x}_-|\le \sigma b_{22},\quad |{\boldsymbol v}_+- R{\boldsymbol v}_-|\le \sigma$.
Here $R(v_x,v_y)=(v_x,-v_y)$ and, as usually, unit vectors ${\boldsymbol v}_{\pm}$ correspond to angles $v_\pm$.

\noindent \textbf{Lemma 2. }\emph{For any $\sigma>0$ there exists a $h_0>0$ such that if $h\in (0,h_0)$ and conditions $(8)$ are satisfied, the Lebesgue measure of the set $N_\sigma$ is greater than $1-\sigma$.}

\noindent\textbf{Proof.} Let $(x_-,v_-)$ be the initial position and the angle of the initial velocity of a particle. We identify $v_-$ with a point of the lower semicircle. Let $(x_+,v_+)$ correspond to the exit of the particle. Here $v_+$ is a point of the upper semicircle.

Observe that any particle, corresponding to initial conditions $(0,v_-)$, $|\tan v_-|>2h$, is reflected back to the same point after two impacts (unless the particle is moving strictly down). Moreover, after the first impact the motion of the particle is strictly parallel to the axis $Ox$. Let $v_0$ and $v_1$ be such that $\tan v_0=-2h$, $\tan v_1=2h$. Then there exists a $\delta>0$ such that every trajectory of the amphora billiard, corresponding to initial conditions $(x,v)$: $|x|<\delta$, $v\in (v_-,v_+)$, $v\neq -\pi/2$ has exactly two impacts and both of them correspond to points of "sides of the amphora"\, i.e. parabolas. It suffices to prove that
$$D_v=\dfrac{\partial (x_+,v_+)}{\partial (x_-,v_-)}(0,v)=\begin{pmatrix} d_{11} & d_{12} \\ d_{21} & d_{22} \end{pmatrix}=
\begin{pmatrix} \pm 1 & 0\\ 0 & -1
\end{pmatrix}$$
for any $v\in (\Theta_-,\Theta_+)$. The sign of the element $d_{11}$ is not important for us.

Since every trajectory that passes via the focus comes back to the focus after two reflections, we have $d_{12}=0$. Due to symmetry reasons, $d_{22}=-1$, $d_{11}=\pm 1$.

Let $\boldsymbol{n}_-$ and $\boldsymbol {n}_+$ be unit normal vectors for points of the first and the second impact respectively. Let $(x_-,v_-)$ be initial conditions for the trajectory. Then ${\boldsymbol n}_\pm$ are functions of $x_-$ and, moreover, grace to the structure of the considered domain, the vector ${\boldsymbol n}_-$ uniquely defines the point of the first impact and, consequently, uniquely defines the vector ${\boldsymbol n}_+$. Let $n_\pm$ be the angles between ${\boldsymbol n}_\pm$ and $Ox$. Consider the angle $\alpha$ between the axis $Ox$ and the trajectory of the particle after the first impact. Clearly, $\alpha=0$ for all solutions, passing via the focus. Due to reflection law, $\alpha=v_--2n_-$. Comparing the trajectory of a particle with one obtained by reversion of time we get $v_--2n_-+v_+-2n_+=\pi$.

On the other hand, for all solutions, passing via the focus, one can easily see that $dn_+/d n_-=-1$. This implies
$$\left.\left(\dfrac{\partial n_+}{\partial x_-}+\dfrac{\partial n_-}{\partial x_-}\right)\right|_{x_-=0}=0 \quad \mbox{and} \quad
\left.\left(\dfrac{\partial v_+}{\partial x_-}+\dfrac{\partial v_-}{\partial x_-}\right)\right|_{x_-=0}=0 \quad \Rightarrow d_{21}=0. \quad \square $$

Note also, that if a trajectory meets the neck of the amphora so that the absolute value of the direction of the entrance velocity is less than $\pi/4$, it is reflected upwards and does not interact with the boundary of the amphora any more.

Amphora billiards have a disadvantage, similar to one of mushrooms: particles can get stuck there, having a big number of impacts until they leave the amphora domain. We modify the amphora in the following way. Attach two triangles $B_LC_LF_L$ and $B_RC_RF_R$ to horizontal parts of the boundary of the billiard (Figure 3(b)). We do it so that
\begin{enumerate}
\item $|B_LF_L|=|B_RF_R|=h^{5/4}\quad (=o(F_LF_R))$,
\item $\angle F_LB_LC_L=\angle F_RB_RC_R=\pi/6$, $\angle F_LB_LC_L=\angle F_RB_RC_R=\pi/4$.
\end{enumerate}

We introduce a coordinate $\Phi\in [-\pi/2,\pi/2]$ on boards of the amphora. This coordinate corresponds to the inclination of the line, passing through the origin and the selected point. Consider two symmetric points $O_L$ and $O_R$ that are centers of segments $[B_LF_L]$ and $[B_RF_R]$ respectively. Replace parts of parabolas, corresponding to $\Phi\in [-\pi/4,\pi/4]$ with arcs of ellipses $E_L$ and $E_R$ such that one focus for both of these ellipses is $O$ and another one is $O_L$ for $E_R$ and $O_R$ for $E_L$. The modified amphora domain is constructed, Figure 3(b).

Now we study billiard trajectories for the modified amphora billiard. Suppose that the angle between the initial velocity and the line is less than $\pi/7$.

If a particle hits the boundary at one of points of $[A_LB_L]$ or $[A_RB_R]$ it is reflected upwards and does not have any other impacts. Otherwise, it interacts twice with arcs of parabolas. After that, due to Lemma 2 there exist following three alternatives, Figure \,3.
\begin{enumerate}
\item A particle leaves the amphora domain forever without having any more impacts.
\item A particle hits $[A_LB_L]$ or $[A_RB_R]$ and leaves the amphora domain.
\item A particle hits $[B_LC_L]$ or $[B_RC_R]$ then maps to a point of $E_L$ or $E_R$ respectively.
\end{enumerate}

Lemma 2 guarantees that the "majority"\ of trajectories behave according to the first scenario. Note that for any initial conditions of the considered type the number of impacts cannot exceed 4.

\section{Hybrid hollows}

\hspace{\parindent} Now we are ready to construct the rough element, i.e. the hollow, corresponding to the rough disc with a prescribed law of reflection. We modify the amphora billiard so that for some selected directions of incident particles it works as a retroreflector and for some others it works as a quasielastic reflector.

\begin{figure}[ht]\begin{center}
\includegraphics*[width=2.5in]{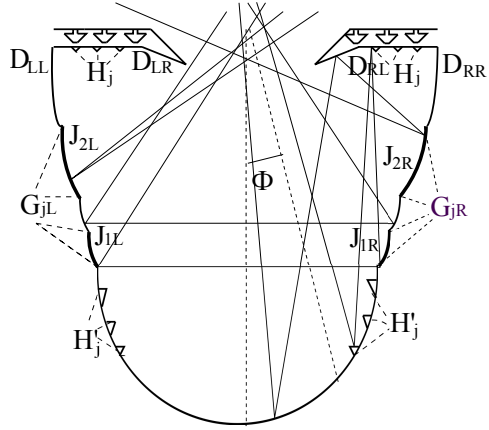}
\end{center}
\caption{Hybrid billiard.}
\end{figure}

Select two symmetric sets of non-intersecting segments $J_{kL}$ and $J_{kR}$ ($k=1,\ldots m$) given by $J_{kR}=[\Phi_k^0,\Phi_k^1]$, $J_{kL}=[-\Phi_k^1,-\Phi_k^0]$. Assume that $5\pi/14=\pi/2-\pi/7<\Phi_1^0<\Phi_1^1<\ldots<\Phi_m^0<\Phi_m^1<\pi/2$.

Given a point $D_{RR}$ (the right edge of the hollow) we attach an arc of the ellipse with the foci at $F_L$ and $F_R$ and corresponding to $\Phi\in (\Phi_m^0,\Phi_m^1)$. Then we draw an arc of the parabola with a focus at the origin and the vertical axis of symmetry through the free end of the constructed arc of the ellipse. We do it for $\Phi\in(\Phi_{m-1}^1,\Phi_m^0)$. We repeat similar constructions of arcs of ellipses with same foci and parabolas with the same focus until we reach $\Phi=\Phi_1^1$. Then we attach the last arc of parabola, corresponding to $\Phi\in (\pi/4,\Phi_1^1)$. To finish the construction we attach an arc of an ellipse, corresponding to $(\pi/4,\pi/2)$ similarly to what we did for modified amphora billiards (Figures\, 3,4).

It may happen that a trajectory or a pseudotrajectory which hits a parabolic part of the boundary near its junction with an elliptic part, next hits an elliptic part on the opposite side of the hollow. Generally, this means that the corresponding billiard trajectory hits one of segments $\Sigma_L=[D_{LL},D_{LR}]$ or $\Sigma_R=[D_{RL},D_{RR}]$ on the upper part of the boundary of the hollow (Figure\,4). Let $G_{1L},\ldots,G_{mL}$ and $G_{1R},\ldots,G_{mR}$ be junction between elliptic and parabolic sectors. Consider $H_1,\ldots H_{2m}$ that are points on the union $\Sigma_L\bigcup\Sigma_R$, corresponding to "parabolic+elliptic"\ reflections from points $G_{kL}$ and $G_{kR}$ or vice versa. We put a system of flat mirrors (segments) of sizes $h^{5/4}$ centered at $H_j$ ($j=1,\ldots,2m$) so that all $h^{3/2}$ pseudotrajectories, hitting first parabolic, then elliptic sectors, are reflected via these mirrors to $h^{9/8}$ neighborhoods of points $H'_1,\ldots H'_{2m}$ such that $H'_j\in (-\pi/4,-\pi/7)\bigcup (\pi/7,\pi/4)$ for all $j$. We put flat mirrors of lengths $h^{9/8}$, centered at points $H'_j$ so that all considered trajectories and pseudotrajectories are reflected by these mirrors to the $h^{17/16}$ neighborhood of the center of the entrance of the hollow (Figure\, 4). That size is still much less than length of the entrance, equal to $h$. Trajectories and pseudotrajectories corresponding to this hybrid billiard, with incident angles $v_-\in (-\pi,-6\pi/7)\bigcup (-\pi/7,0)$ are the following.
\begin{enumerate}
\item If a pseudotrajectory does not hit points, corresponding to one of segments $J_{kL}$ or $J_{kR}$, the behavior is the same as for the modified amphora billiard.
\item If it hits one of the mentioned segments, is reflected "almost back"\ (similarly to what happens for Bunimovich mushrooms). Then the \emph{pseudotrajectory} leaves the domain without farther interactions with walls.
\item A small proportion of particles (which tends to 0 as $h\to 0$) has a distinct behavior. However, all such particles leave the hollow, having at most 4 impacts.
\end{enumerate}

So, the constructed hollow is perfect. Now we describe how it is possible to cover almost all segment $I\in Ox$ (we may also do the same if $I$ is an arc of the circle) with tops of hybrid billiards. Cut the middle part of $I$ of the length $2 h |I|$ and insert there a hybrid billiard of imperfectness $h$ and the basis of the neck equal to $2h |I|$. Call this hollow one of the first generation. Let $b_1$ be the corresponding rescaling coefficient. Take $b_2=\varrho h^2 b_1$. Here $\varrho<1$ is the principle rescaling for smaller mushrooms of the "second generation", Figure 4). Then we put $N\sim {h}^{-1}|I|$ non-intersecting hollows of the second generation whose tops correspond to subsegments of $I$. We repeat this procedure, creating hollows of the third level and so on. On the step number $L$, the measure of the part of the segment $I$, not covered by tops of already constructed hollows can be estimated by the value
$|I|(1-{\widetilde h}/2)^L$. In the limit, we get a zero-measure Cantor set. However, we stop after finitely many steps.

\section{Main result}

\noindent\textbf{Theorem 1.} \emph{Let $\boldsymbol{g}: [a,\, b] \to {\mathbb{R}}^2$ be a continuous rectifiable curve. Then for any $\varsigma > 0$ there exists a motion $(\boldsymbol{X}(\tau),\varphi(\tau)$ of a rough disc with a radius $r > 0$ and with the coordinate of center $\boldsymbol{X}(\tau)$ such that after a continuous and monotone increasing change of parameter $\tau=\tau(t)$,\, $t \in [a,\, b]$ one has
$$|{\boldsymbol g}(t) - \boldsymbol{X}(\tau(t))| < \varsigma.\eqno (9)$$
Here $\boldsymbol{X}(\tau)$ is the position of the center of the disc; $\varphi(\tau)$ is the turn of the disc.}

Note that we the curve ${\boldsymbol g}$ is not necessarily injective: self-intersections and even coincidence of some fragments of the curve are allowed.

The following auxiliary theorem states that any broken line can be approximated by trajectories of rough discs. Namely, let ${\boldsymbol G}(t)$,\, $t \in [a,\, b]$ be a parameterized broken line with a finite number of segments, $\Gamma=\{\boldsymbol{G}(t):t\in[a,b]\}$.  Self-intersections are allowed, but we require that no vertex of the broken line is a point of intersection. Moreover, we approximate broken lines so that inclinations of every segment with respect to the previous one varies from $-\pi/4$ to $0$. For instance, instead a rotation by the angle $\pi/4$, we apply seven rotations by $-\pi/4$. 

\noindent\textbf{Theorem 2.} \emph{For any $\varsigma > 0$ there exists a motion of a rough disc of radius $r > 0$ whose center is $\boldsymbol{X}(\tau)$ such that after a continuous and monotone increasing change of parameter $\tau=\tau(t)$,\, $t \in [a,\, b]$ inequality $(9)$ is satisfied.}

Theorem 1 is an obvious consequence of Theorem 2. Indeed, each rectifiable curve can be uniformly approximated by broken lines, and each broken line can be uniformly approximated by trajectories of rough discs. 

\noindent\textbf{Proof of Theorem 2}.

First we notice that a curve homothetic to the trajectory of a rough disc is also a trajectory of a rough disc. Let $\boldsymbol{X}(t)$ be the motion of the center of a rough disc of radius $r$. Let $\omega(t)$ be its angular velocity and $\epsilon$ be a positive constant. Then the coordinate of the center of a disc of radius $\nu r$ homothetic to the original one moving in the same medium with the initial velocity $\epsilon \boldsymbol{X}'(0)$ and the initial angular velocity 
$\omega(0)$, is given by $\epsilon\boldsymbol{X}(t)$, and its angular velocity is $\omega(t)$.

This scaling argument allows one to reduce Theorem 2 to the problem of approximation of a broken lines $\frac{1}{\epsilon} \boldsymbol{g}(t)$ where $\epsilon$ is a small parameter. Select a splitting of the broken line into segments with ends, corresponding to $a=T_0<T_1<\ldots<T_{m-1}<T_m=b$.

Take a disc $B_{n_\varepsilon}$ with the roughness of the considered form. Introduce the measure in $[-\pi/2,\, \pi/2] \times [-\pi/2,\, \pi/2]$ associated with the cavity which has the density
$$\frac 12 \cos x\{ \delta(x-y) \cdot \chi_{J\cup J'}(x) + \delta(x+y) \cdot \big[ 1 - \chi_{J\cup J'}(x) \big]\}\, dx\, dy. \eqno (10)$$
if $|x|,|y|\le 5\pi/14$. Here
$$J = \bigcup\limits_{i=1}^m J_{iR}=
\bigcup\limits_{i=1}^m \big[\pi/2 - e^{-T_i/\varepsilon},\, \pi/2 - e^{-(T_i+\Delta T_i)/\varepsilon}\big] \quad \text{and} \quad J' = -J,\eqno (11)$$
$i = 1,\, 2,\ldots,m$ is a finite set of indices. We select $J_{iR}$ as "elliptic"\ segments on the boundary of a cavity (see Section 7 and Figure 4).
The initial angular velocity $\lambda(0) = \omega(0)/rV(0)$ is taken to be $\lambda(0) = e^{T_0/\varepsilon}$. $T_i - T_{i-1}$ is the length of the $i$~--~th segment of the broken line, $\Delta T_i = \varphi_i e^{-T_i/\varepsilon}$. Here $\varphi_i$ are parameters, close to angles $\varphi^0_i\in [-\pi/4,0]$ between the $i$-th and $(i+1)$-th segments of the broken line. Now we note that $\varepsilon>0$ is taken so small that all segments $J_i$ are disjoint.

As the disc moves, the relative angular velocity increases and less of the part of the cavity is "observable"\ by particles. Depending on the value of $\lambda$ either $J\bigcup J'$ or completion of this set dominate in the "observable"\ part. Respectively, we have rotation or "almost straight forward"\ motion. A small part $\varepsilon$ of the boundary is filled with cavities. The rest, $1-\varepsilon$, of the boundary is not filled, that is, is just a union of arcs of the unit circumference. Both parts are uniformly distributed along the boundary.

Consider the natural parametrization ${\boldsymbol g}(\tau)$,\, $\tau \in [T_0,\, T_m]$, where $[T_{j-1},\, T_j]$ parameterize segments of the broken line. We find a motion of a rough disc of unit radius where $\boldsymbol{X}(\tau)$ is the position of the center ,\, and values $S_j$ ($j=0,\ldots,m$),
$\tau \in [S_0,\, S_m]$ such that
$$|{\boldsymbol g}(\tau)/\varepsilon -{\boldsymbol X}(\tau/\varepsilon)| < (m+1)/\sqrt\varepsilon\eqno (12)$$
or, equivalently, $\left|{\boldsymbol g}(\tau) -\varepsilon{\boldsymbol X}(\tau/\varepsilon)\right| < (m+1)\sqrt\varepsilon$. Given $\varsigma$, we take $\varepsilon$ so that $\varsigma=(m+1)\sqrt{\varepsilon}$. Then we take the rescaling parameter $\kappa=\varepsilon$ and easily obtain inequality $(9)$.

The motion of the disc is described in terms of the parameter $\tau$ proportional to the natural one (see $(5)$). It can be deduced from equations $(6)$ and $(12)$ and from equations defining the measure $(10)$, $(11)$ that the differential equation for $\lambda(\tau)$ takes the form $\lambda' = \lambda u(\lambda,\varepsilon,\tau)$ where $u(\lambda,\varepsilon,\tau)\rightrightarrows 1$ as $\varepsilon \to 0$. So $\lambda = e^{w(\varepsilon,\tau)}$ where $w$ is increasing with respect to $\tau$ and $w(\varepsilon,\tau)/\tau\rightrightarrows 1$. Consider values $S_j$ defined by equalities $w(\varepsilon,S_j)=T_j/\varepsilon$.

Using equations $(3)$, $(4)$ and $(6)$, introduce the notation $x_0 = x_0(\lambda) = \arccos(1/\lambda)$, and obtain the equality
$$\frac 12\int_{x_0}^{\pi/2} c_T(x,-x,\lambda)\, \cos x\, dx = 0$$
(recall that the function $c_T$ is defined by $(4)$). This means that the component, orthogonal to the current velocity, of the force acting on a smooth (without roughness) disc is zero. So we obtain $\theta'(\tau) = -\varepsilon R_T(\lambda(\tau))$
where
$$R_T(\lambda) = \frac 12\int_{[x_0,\pi/2]\cap J} (c_T(x,x,\lambda) - c_T(x,-x,\lambda)) \cos x\, dx,\eqno (13)$$
with $c_T(x,x,\lambda) - c_T(x,-x,\lambda) =$
$$\frac{3\sin x}{\sin\zeta} \{ (\lambda^3\sin^3x + 3\lambda\sin x\sin^2\zeta) \cos\zeta\sin x + (3\lambda^2\sin^2x\sin\zeta + \sin^3\zeta) \sin\zeta\cos x \} $$
and $\zeta = \zeta(x) = \arccos(\lambda\cos x)$. After some algebra we get
$$c_T(x,x,\lambda) - c_T(x,-x,\lambda)
= \frac{3\sin x\cos\zeta}{\lambda\sin\zeta} \{ (\lambda^2 - \cos^2\zeta)^2 + 6\sin^2\zeta(\lambda^2 - \cos^2\zeta) + \sin^4\zeta \}.$$

Making the change of variable $x \to \zeta$ in the integral $(13)$, we obtain
$$R_T(\lambda) = \int_{[0,\pi/2]\cap \tilde J} \frac{3}{2\lambda^3}\, \{ (\lambda^2 - \cos^2\zeta)^2 + 6\sin^2\zeta(\lambda^2 - \cos^2\zeta) + \sin^4\zeta \}\, \cos^2\zeta\, d\zeta,\eqno (14)$$
where
$$\tilde J = \bigcup_{j=0}^{m-1}[\zeta_j,\, \zeta_j + \Delta\zeta_j],$$
with $\zeta_j = \arccos(\lambda e^{-S_j})$,\, $\zeta_j + \Delta\zeta_j = \arccos(\lambda e^{-w^{-1}((T_j+\Delta_j)/\varepsilon})$. Notice that the expression $\{\ldots\}$ in the integral in the right hand side of $(14)$ can be estimated as $\{\ldots\} = \lambda^4+O(\lambda^3)$ for large values of $\lambda$.

Substituting $\lambda = e^{w(\varepsilon,\tau)}$, one obtains
$$\zeta_j = \arccos(e^{w(\varepsilon,\tau)-T_j/\varepsilon}) \quad \text{and} \quad \Delta\zeta_j = \frac{e^{w(\varepsilon,\tau)-T_j/\varepsilon}}{\sqrt{1 - e^{2w(\varepsilon,\tau)-2T_j/\varepsilon}}}\, \frac{\Delta_j}{\varepsilon}(1+o_\varepsilon(1))$$
where $o_\varepsilon(1)\to 0$ as $\varepsilon\to 0$. The value of $\varepsilon R_T(\lambda)$ can now be evaluated as
$$\varepsilon R_T(\lambda) = \varepsilon\, \frac{3}{2\lambda^3}\, (\lambda^4 \cos^2\zeta_j \Delta\zeta_j +{\hat R}^0_j(\lambda,\varepsilon)) =$$
$$\varepsilon\, \frac{3\lambda}{2}\, e^{2w(\varepsilon,\tau)-2T_j/\varepsilon}\, \frac{e^{w(\varepsilon,\tau)-T_j/\varepsilon}}{\sqrt{1 - e^{2w(\varepsilon,\tau)-2T_j/\varepsilon}}}\, \frac{\Delta_j}{\varepsilon} +{\hat R}^1_j(\tau,\varepsilon)= \frac{3\varphi_j}{2}\, \frac{e^{4w(\varepsilon,\tau)-4T_j/\varepsilon}}{\sqrt{1 - e^{2w(\varepsilon,\tau)-2T_j/\varepsilon}}}+{\hat R}^1_j(\tau,\varepsilon).$$
Here $|{\hat R}^0_j(\lambda,\varepsilon)|\le C\lambda^3$ where $C$ is a constant; ${\hat R}^1_j(\tau,\varepsilon)$ tends to zero as $\lambda(\tau)\to\infty$, $\varepsilon\to 0$. Thus, we come to the following differential equation for $\theta(\tau)$,
$$\frac{d\theta}{d\tau} =\dfrac{3\varphi_j}{2}\, \dfrac{e^{4w(\varepsilon,\tau)-4T_j/\varepsilon}}{\sqrt{1 - e^{2w(\varepsilon,\tau)-2T_j/\varepsilon}}}+
{\hat R}^1_j(\tau,\varepsilon), \qquad \text{if} \quad \tau \in [S_j,\, S_{j+1}-1/\sqrt{\varepsilon}].$$

Solutions for this equation are
$$\theta(\tau) = \theta(S_j) + \varphi_j \Big[ 1 - \sqrt{1 - e^{2w(\varepsilon,\tau)-2T_j/\varepsilon}} \Big( 1 + \frac 12\, e^{2w(\varepsilon,\tau)-2T_j/\varepsilon} \Big)\Big]+{\widetilde R}(\varepsilon,\tau),
\eqno (15)$$
$$\text{if} \quad \tau \in [S_j,\, S_{j+1}];\quad j=0,\ldots,m-1.$$

Here $|\widetilde R(\varepsilon,\tau)|\le \sqrt{\varepsilon}$ if $\varepsilon$ is sufficiently small. The function $\theta$ is increasing with respect to $\tau$ and with respect to each parameter $\varphi_j$. So, we can select all $\varphi_j$ so that $\theta(S_{j+1})-\theta(S_j)=\varphi^0_j$. Thus, any part of the trajectory ${\boldsymbol X}([S_j,\, S_{j+1}-1/\sqrt{\varepsilon}])$ ($j\ge 1$) is an arc, close to a line segment of length $(T_{j+1} - T_{j})\varepsilon^{-1}-\varepsilon^{-1/2}$.

Let $L$ be the length of the curve ${\boldsymbol g}$, $\theta_0(\tau)$ be the piecewise constant function, equal to $0$ on $[S_0,S_1)$ and equal to
$\varphi_1^0+\ldots+\varphi_{j-1}^0$ on $[S_{j-1},S_j)$. Then for any $\tau\in [0,L]$ we have
$$\left|{\boldsymbol X}(\tau/\varepsilon) - \frac{1}{\varepsilon} {\boldsymbol g}(\tau)\right|\le \int_{S_0}^{\tau/\varepsilon} |\varphi(s)-\varphi_0(s)|\,ds.
\eqno (16)$$

For $\tau \in [\varepsilon S_{j},\, \varepsilon S_{j+1}-\sqrt\varepsilon]$ the velocity vector ${\boldsymbol X}'(t/\varepsilon)$ forms an angle $O(e^{-\frac{1}{2\sqrt\varepsilon}})$ with the $j$~--th segment of the broken line. This follows from representations $(15)$. On the other hand, contributions of any segment $[\varepsilon S_{j}-\sqrt\varepsilon,\varepsilon S_{j}]$ to the right hand side of $(16)$ are estimated by $\varepsilon^{-1/2}$. So,
we have inequality $(12)$ satisfied if $\varepsilon$ is small. $\square$

\section{Conclusion and discussion}

\hspace{\parindent} The main results of this paper are the following. Two-dimensional trajectories of bodies, whose boundaries are close to circles, may have (up to rescaling) any shape. The same statement is true for flat curves in the three dimensional real space.  Also a description of amphora billiard (quasi-elastic reflector) and its modifications with a wide variety of response functions have been given. All these results are principally novel.

However, our construction while being mathematically correct cannot be implemented in practice. First, we make some non-realistic assumptions that the medium temperature is absolute zero, the particles of the medium do not collide, and (even worse) the collisions of the particles with the boundary of the body are perfectly elastic. Second, even if all these assumptions are satisfied, each cavity should be fabricated with exceptionally high precision, the scale of precision being much smaller than the size of atoms. Third, the path traversed by a disc is proportional to the logarithm of time. Roughly speaking, it may happen that the first meter of the trajectory is traversed in a second, the second meter in a minute, the third meter in a hour, ..., the tenth meter in a billion of years. The experimenter may just not survive the end of the experiment.

Imagine a football player who wants to send the ball so that the trajectory goes round all the players of the rival team and finally gets into the gate. He can indeed do so making use of our results, but the ball surface should be very special; the pressure of the atmosphere should be very low; the Earth gravitation should be negligible; the rival players should be asked not to prevent the (eventually very small) motion of the ball. And it remains to wait. Oh, forgot to say that all this should happen in two dimensions.

\end{document}